\documentclass{imsart}
\RequirePackage[OT1]{fontenc}
\RequirePackage{amsthm,amsmath,amssymb}
\RequirePackage[numbers]{natbib}
\RequirePackage[colorlinks,citecolor=blue,urlcolor=blue]{hyperref}

\arxiv{arXiv:0000.0000}

\startlocaldefs
\numberwithin{equation}{section}
\theoremstyle{plain}
\newtheorem{Theorem}{{\sc THEOREM}}[section]  
\newtheorem{Lemma}[Theorem]{{\sc LEMMA}}

\newtheorem{Definition}[Theorem]{{\sc Definition}}

\endlocaldefs

\def\endpf{\hfill $\Box$}

\newcommand {\bbeta}{\mbox{\boldmath $\beta$}}
\newcommand {\bbR}{{\mathbb R}}

\newcommand {\bbN}{{\mathbb N}}

\newcommand {\bbZ}{{\mathbb Z}}
\newcommand {\bbC}{{\mathbb C}}

\newcommand {\iy}{\infty}
\newcommand {\FF}{\mathcal{F}}

\newcommand {\bbJ}{\1}

\newcommand \mes{{\rm \lambda}}

\def\bt{{\bf t}}
\def\bs{{\bf s}}
\def\bm{{\bf m}}
\def\bx{{\bf x}}

\def\eps{{\varepsilon}}
\def\-{\!-\!}
\def\={\!=\!}
\def\+{\!+\!}

\def\1{\hbox{\upshape1\kern-.15em\vrule height 1.6ex width .3pt%
\vrule width .8pt height .25pt\kern.15em}}
\def\mod{~ {\textsf{\upshape mod}}~1} 

\def\modb{~ {\textsf{\upshape mod}}~{\bf 1}}

\begin{document}

\begin{frontmatter}
\title{Equidistribution, Uniform distribution: a probabilist's perspective}
\runtitle{Equidistribution, Uniform distribution: a probabilist's perspective}

\begin{aug}
\author{\fnms{Vlada} \snm{Limic}\ead[label=e1]{vlada@math.unistra.fr}}
\and
\author{\fnms{Ned\v{z}ad} \snm{Limi\'{c}}\thanksref{t3}\ead[label=e2]{nlimic@math.hr}}

\address{IRMA, UMR 7501 de l’Universit{\'e} \\
de Strasbourg et du CNRS,\\
7 rue Ren{\'e} Descartes,\\67084 Strasbourg Cedex, France \\
\printead{e1}}

\address{University of Zagreb,\\
Bijeni\v{c}ka cesta 30,\\10000 Zagreb, Croatia \\
\printead{e2}}

\thankstext{t3}{supported in part by the Croatian Science Foundation under project 9780 WeConMApp}
\runauthor{V. Limic and N. Limi\'{c}}

\affiliation{Some University and Another University}

\end{aug}

\begin{abstract}
The theory of equidistribution is about hundred years old, and has been developed primarily by number theorists and theoretical computer scientists. 
A motivated uninitiated peer could encounter difficulties perusing the literature, due to various synonyms and polysemes used by different schools.
One purpose of this note is to provide
a short introduction for probabilists.  
We proceed by recalling a perspective originating in a work of the second author from 2002.
Using it, various new examples of completely uniformly distributed $\!\!\mod$ sequences, in the ``metric'' (meaning almost sure stochastic) sense, can be easily exhibited. In particular, we point out natural generalizations of the original $p$-multiply equidistributed sequence $k^p\, t \mod$, $k\geq 1$ (where $p\in \bbN$ and $t\in[0,1]$), due to Hermann Weyl in 1916. 
In passing, we also derive a Weyl-like criterion for weakly completely equidistributed (also known as WCUD) sequences, of substantial recent interest in MCMC simulations.

The translation from number theory to probability language brings into focus a version of the strong law of large numbers for weakly correlated complex-valued random variables, the study of which was initiated by Weyl in the aforementioned manuscript, followed up by Davenport, Erd\H{o}s and LeVeque in 1963, and greatly extended by Russell Lyons in 1988.
In this context, 
an application to $\infty$-distributed Koksma's numbers $t^k \mod$, $k\geq 1$ (where $t\in[1,a]$ for some $a>1$), and an important generalization by Niederreiter and Tichy from 1985 are discussed.

The paper contains negligible amount of new mathematics in the strict sense, but its perspective and open questions included in the end could be of considerable interest to probabilists and statisticians, as well as
certain computer scientists and number theorists.
\end{abstract}
\end{frontmatter}

\section{Introduction}\label{sec1}
This is certainly neither the first nor the last time that equidistribution is viewed using a ``probabilistic lense". 
A probability and number theory enthusiast will likely recall in this context the famous Erd\"os-Kac central limit theorem \cite{erdkac} (see also \cite{durrett}, Ch.~2 (4.9)), or the celebrated monograph by Kac \cite{kac}.

Unlike in \cite{kac,kac_paper,erdkac,kesten1,kesten2} our main concern here are the completely equidistributed sequences and their ``generation''. 
In the abstract we deliberately alternate between equidistribution (or equidistributed) and two of its synonyms. 
In particular, the {\em completely uniformly
 distributed} (sometimes followed by $\!\!\mod$) and {\em $\infty$-distributed} in the abstract, the keywords, and the references 
(see for example any of
 \cite{hol1,hol2,drmota_tichy,knu_paper,knu69,kor92,kuni,nietic,lacaze,loynes,levinM,tribble_owen,owen_tribble,strpor})
mean the same as {\em completely equidistributed}.
Uniform distribution is a fundamental probability theory concept. To minimize the confusion, in the rest of this note  we shall: \\
(i) always write the {\em uniform law} when referring to the distribution of a uniform random variable, and \\
(ii)  almost exclusively write equidistributed (to mean equidistributed, uniformly distributed,  uniformly distributed $\!\!\mod$ or $\,\cdot\,$-distributed), usually preceded by one of the following attributes: {\em simply, $d$-multiply or completely}. 

\smallskip
Let $D = [0,1]^d$ be the $d$-dimensional unit cube. 
Let $r\in \bbN$ and assume that a bounded domain $G$ is given in ${\bbR}^r$.
Denote by $\mes$ the restriction of the $d$-dimensional Lebesgue measure on $D$, as well as the Lebesgue measure on $\bbR^r$. In particular $\mes(G)$ is the $r$-dimensional volume of $G$. In most of our examples  $r$ will equal $1$, and $G$ will be an interval. 

Throughout this note, a sequence of measurable (typically continuous)  functions $({\bf x}_k)_{k\geq 1}$, where ${\bf x}_k:G \mapsto {\bbR}^d$  will define a sequence of points $\mbox{\boldmath $\beta$}_k \in D$ as follows:
\begin{equation}
\label{D:beta}
 \mbox{\boldmath $\beta$}_k \equiv  \bbeta(\bt)  :=  {\bf x}_k({\bf t}) \modb, \quad
 {\bf t} \in G, 
\end{equation}
understanding that ${\bf 1}=(1,\ldots,1)\in \bbR^d$, and that the modulo operation is naturally extended to vectors in the component-wise sense. 
While each $\bbeta_k$ is a (measurable) function of $\bt$,  this dependence is typically omitted from the notation.
The $({\bf x}_k)_k$ is called the {\em generating sequence} or simply the {\em generator}, and the elements $\bt$ of $G$ are called {\em seeds}.

Any $A = \prod_j (a_j, b_j] \subset D$, $ 0 \leq
a_j < b_j \leq 1$ will be called a {\em $d$-dimensional box}, or just a {\em box} in $D$. 
As usual, we denote by ${\bbJ}_S$ the indicator of a set $S$.
The sequence $\mbox{\boldmath $\beta$} =
\{\mbox{\boldmath $\beta$}_k : k \in {\bbN}\}$ is {\em simply equidistributed in $D$} if
\begin{equation}
\label{EeqWeyl} \ \lim_{N \to \infty}\:
 \frac{1}{N} \sum_{k = 1}^N \:{\bbJ_A} (\mbox{\boldmath $\beta$}_k(\bt)) \ = \ 
 \mes(A),
\end{equation}
for each $A$ box in $D$, and almost every $\bt\in G$.

\subsection{Notes on the literature}
\label{S:lite}
Our general setting is mostly inherited from \cite{comed}. 
The few changes in the notation and the jargon aim to simplify reading of the present work by an interested probabilist or statistician. 
In \cite{drmota_tichy,knu69,kuni,kor92,strpor} and other standard references, equidistribution is defined as a property of a single (deterministic) sequence of real numbers. 
Nevertheless, in any concrete example discussed here (e.g.~the sequences mentioned in the abstract, as well as all the examples given in the forthcoming sections), this property is verified only up to a null set over a certain parameter space. So there seems to be no loss of generality in 
integrating the almost everywhere/surely aspect in Definitions \ref{def1} and \ref{def2} below.
There is one potential advantage: the Weyl criterion (viewed in a.s.~sense) reduces to (countably many applications of) the strong law of large numbers (SLLN) for specially chosen sequences of dependent complex-valued random variables. The just made observation is the central theme of this note.  
Note in addition that the study of R4 (and R6) types of randomness (according to Knuth \cite{knu69}, Section 3.5, see also Sections \ref{S:discNT} and \ref{S:conclOP} below) makes sense only in the stochastic setting.

\smallskip
Davenport, Erd\"os and LeVeque \cite{delv63} are strongly motivated by the Weyl equidistribution analysis \cite{We}, yet they do not mention any connections to probability theory.
Lyons \cite{lyo88} clearly refers to the main result of \cite{delv63} as a SLLN criterion, but is not otherwise interested in equidistribution. 

The Weyl variant of the SLLN (see Section \ref{S:weylSLLN}) is central to the analysis of Koksma \cite{koksma}, who seems to ignore its probabilistic aspect.
The breakthrough on the complete equidistribution of Koksma's numbers (and variations) by Franklin \cite{F2} relies on the main result in \cite{koksma}, but without (any need of) recalling the Weyl variant of the SLLN in the background (see the proof of \cite{F2}, Theorems 14 and 15). On the surface, this line of research looks more and more distant from the probability theory.
The classical mainstream number theory textbooks (e.g.~\cite{kuni}), as well as modern references (e.g.~\cite{drmota_tichy,strpor}), corroborate this point of view of  equidistribution in the ``metric'' (soon to be called ``stochastic a.s.'') sense.
And yet, the Niederreiter and Tichy \cite{nietic} metric theorem, considered by many as one of the highlights on this topic, consists of lengthy and clever (calculus based) covariance calculations,  followed by an application of the SLLN from \cite{delv63} (see Section \ref{S:discNT} for more details).  

To the best of our knowledge, \cite{comed} is the first (and at present the only) study of complete (or other) equidistribution in the ``metric'' sense, which identifies the verification of Weyl's criterion as a stochastic problem (equivalent to countably many SLLN, recalled in Section \ref{S:equi_prob}), without any additional restriction on the nature of the generator $({\bf x}_k)_k$.
In comparison,  Holewijn \cite{hol1,hol2} made analogous connection (and in \cite{hol2} even applied the SLLN criterion of \cite{delv63} to the sequence of rescaled Weyl's sums (\ref{EidenW})), but only under particularly nice probabilistic assumptions, which are not satisfied in any of the examples discussed in the present survey.
In addition, both Lacaze \cite{lacaze} and Loynes \cite{loynes} 
use the Weyl variant (or a slight modification) of the SLLN, again under rather restrictive probabilistic hypothesis (to be recalled in Section \ref{S:why_not}).

Kesten in \cite{kesten1,kesten2}, as well as in his subsequent articles on similar  probabilistic number theory topics, works in the setting analogous to that of \cite{hol1,hol2}, however his analysis is not directly connected to the Weyl criterion.
On a related topic, in a pioneer study of (non-)normal numbers and their relation to equidistribution, Mendes France \cite{MF}  also applies the SLLN from \cite{delv63} as a purely analytic result, though it is not clear that a probabilistic interpretation brings a new insight in that setting.
Interested readers are referred to \cite{kemperman,khoshnevisan} for probabilistic perspectives on normal numbers.

\smallskip
The study of complete equidistribution in the standard (deterministic) sense, and its close link to normal numbers, 
was initiated by Korobov \cite{kor48} in 1948. 
Korobov used Weyl's criterion \cite{We} for (multiple and complete) equidistribution, and exhibited an explicit function $f$ (in a form of a power series) such that $(f(k) \mod)_k$ is a completely equidistributed sequence (see also \cite{kor92}, Theorem 28).
Knuth was aware of \cite{F2}, but apparently unaware of either \cite{We,kor48}, when he exhibited in \cite{knu_paper} a (deterministic) completely equidistributed (there called ``random'') sequence of numbers in $[0,1]$, by extending the method of Champernowne's that previously served to find an explicit normal number. 
Knuth \cite{knu_paper,knu69} uses his own (computer science inspired) criterion for complete equidistribution.
All the examples given in \cite{kor48,kor49} and \cite{knu_paper,knu69}, as well as those obtained later on by the ``Korobov school''  (see Levin \cite{levinM}, Korobov \cite{kor92}, and also the historical notes in \cite{kuni}, and \cite{strpor}), are practical to a varying degree. To find out more about the deterministic setting, readers are encouraged to use the pointers (to synonyms and references) given above as a guide to the literature. 

\smallskip
We finally wish to make note of a recent spur of interest in completely equidistributed sequences, and their generalizations (definable only in the stochastic setting), in relation to Markov chain Monte-Carlo simulations \cite{chen_etal,chen_dick_owen,owen_tribble,tribble_owen}. 
Section \ref{S:WCUD} makes a brief digression in this direction.
Importance of complete equidistribution for MCMC is not surprising, in view of a long list of empirical tests that these sequences satisfy (see \cite{F2,knu69} and Remark 2(d,e)). \\
The first author is convinced that anyone who regularly runs or even looks at pseudo-random simulations should benefit from reading a note of this kind. 

\medskip
{\bf Disclaimer and motivation} This review and tutorial is highly inclusive, but by no means exhaustive. The theory of equidistribution (or uniform distribution) is rich, complex and fast evolving, and it would be very difficult to point to a single book volume, let alone a survey paper, which covers all of its interesting aspects (for example, 
the list of references in the specialized survey \cite{aisber_review} overlaps with ours in only four items). Even when focusing on complete equidistribution in the metric (stochastic a.s.) sense, it seems hard to find a single expository article aimed at specialists, let alone at the probability and statistics community at large. We hope to have accomplished here an ``order of magnitude'' effect, citing several dozens of original research papers and surveys, textbooks or monographs, in a brief attempt to shed a probability-friendly light on the concepts and ideas presented, as well as to point out a number of natural and interesting open 
questions. We wish we had come across such a paper at our very encounter with this important topic.
The latter thought gave the impetus to our writing.

\subsection{One-dimensional examples}
Suppose that $d=r=1$, and that $G$ is an interval.
We recall several well-known examples of functions
$x_k: G\mapsto {\bbR}$, that generate equidistributed sequences in $[0,1]$.
The first class of examples is as follows: for each  $p \in {\bbN}$ define
\begin{equation}
\label{ex1.1}
 x_k(t) \ = \ k^p\,t, \quad t \in G \::=\:(0,1),\ k\geq 1.
\end{equation}
Then $(\beta_k)_{k\geq 1} \in [0,1]^\bbN$ from (\ref{D:beta}) are {\em Weyl's numbers}  of parameter $p$.
The  second class of examples are the so-called  {\em multiplicatively generated numbers}: for each $M \in \{2,3,\ldots,\}$ let
\begin{equation}
\label{ex1.2}
 x_k(t) \ = \ M^k \,t, \quad t \in G \::=\:(0,1), \ \ k\geq 1,
\end{equation}
and $(\beta_k)_{k\geq 1} \in [0,1]^\bbN$ as in (\ref{D:beta}).
Finally, if one defines for each $k$, and any $a > 1$
\begin{equation}
\label{ex1.3}
 x_k(t) \ = \ t^k, \quad t \in G \::=\:(1,a), \ \ k\geq 1,
\end{equation}
then $(\beta_k)_{k\geq 1} $ from  (\ref{D:beta}) are known as {\em Koksma's numbers}. 
It is well-known that for each of the above three examples, the sequences $(\beta_k)_k$ are (at least simply)
equidistributed for almost all seeds \cite{We,F2,koksma} (see also \cite{kuni,drmota_tichy,strpor}).
In fact if $p=1$, then the Weyl (Sierpinski, Bohl) equidistribution theorem is a stronger claim: $k\cdot t \mod$ is simply equidistributed for all irrational $t$.

One could {\bf refine} the notion of a {\bf seed}, and call $t$ a seed only if $\beta_k(t)$ is (sufficiently) equidistributed. The new set of seeds $T$ would then be well-defined up to a null-set.
Note that (\ref{ex1.1}) and (\ref{ex1.2}) are both linear in $t$, 
which is not true for (\ref{ex1.3}). 
The methodology developed by the second author in \cite{comed} was motivated by the particularly simple analysis of the linear case, which can be extended (under certain hypotheses on $(x_k)_k$) to the non-linear setting.

\subsection{Multiple equidistribution with examples}
Let $d\geq 2$ be fixed.
A sequence of real measurable functions $(x_k)_{k\geq 1}$ can be used to form sequences $({\bf x}_k)_{k\geq 1}$ of $d$-dimensional vector-valued (measurable) functions, and therefore the corresponding  sequences
$(\bbeta_k)_{k\geq 1}$  (see (\ref{D:beta})), in at least three different natural ways:
\begin{itemize}
\item[{(a)}]
The set of seeds could be equally $d$-dimensional (here $r=d$). More precisely, for ${\bf t} = (t_1, t_2,\ldots, t_d)\in G$, let 
$
{\bx}_k(\bt):= ( x_{k1} (\bt),x_{k2}(\bt),\ldots, x_{kd}(\bt)),
$
where $x_{kj}(\bt) := x_{(k-1)d+j}(t_j)$, for all $k\in \bbN$, $j=1,\ldots, d$. In this case 
\begin{equation}\label{ex1.4mult}
 \beta_{kj}^d \ = \ x_{(k-1)d+j}(t_j) \mod, \quad k\in \bbN, \,j=1,\ldots, d .
\end{equation}
\item[{(b)}]
The set of seeds could be $r$-dimensional, and the successive $\bx$ (and $\bbeta$) could be formed by shifting the window of observation by $1$. More precisely, for $\bt \in G \subset \bbR^r$, define 
$
{\bx}_k(\bt):= ( x_{k1} (\bt),x_{k2}(\bt),\ldots, x_{kd}(\bt)),
$
where $x_{kj}(\bt) := x_{k+j-1}(\bt)$, for all $k\in \bbN$, $j=1,\ldots, d$. In this case 
\begin{equation}\label{ex1.5mult}
 \beta_{kj}^d \ = \ x_{k+j-1}(\bt) \mod, \quad k\in \bbN,\, j=1,\ldots, d.
\end{equation}
\item[{(c)}]
Let $d\geq 1$. The set of seeds could be $r$-dimensional, and the successive $\bx$ (and $\bbeta$) could be formed by shifting the window of observation by $h\in \bbN$. More precisely, for $\bt \in G \subset \bbR^r$, define 
$
{\bx}_k(\bt):= ( x_{k1} (\bt),x_{k2}(\bt),\ldots, x_{kd}(\bt)),
$
where $x_{kj}(\bt) := x_{(k-1)h+j}(\bt)$, for all $k\in \bbN$, $j=1,\ldots, d$. In this case 
\begin{equation}\label{ex1.6mult}
 \beta_{kj}^{d,h} \ = \ x_{(k-1)h+j}(\bt) \mod, \quad k\in \bbN,\, j=1,\ldots, d.
\end{equation}
\end{itemize}
Note that class (c) comprises class (b) (if $h$ is set to $1$), and that if $h>d$,
 then $(\bbeta_k^{d,h})_k$ is formed from a strict subsequence of $(x_k)_k$.
We shall consider the above definitions with $d$ varying over $\bbN$.
The sequences in (\ref{ex1.4mult}) will be included in the analysis of Section \ref{S:linear}.

The equidistribution analysis is typically done on the sequences of vectors $\bbeta^d$ from (\ref{ex1.5mult}). 
Yet (see Remark 1(c) below) the construction in (\ref{ex1.6mult}) is particularly interesting from the perspective of comparison with (pseudo-)random simulations.  

Simple equidistribution in $D$ was defined in the paragraph preceding Section \ref{S:lite}.
Let again $G$ be a bounded domain in $\bbR^r$ for some $r\geq 1$.
\begin{Definition}
\label{def1}
Assume that $(x_k)_{k\geq 1}$ is a sequence of real measurable functions on $G$. 
(i) The sequence $\{\beta_k : k \in {\bbN}\}$ defined in (\ref{D:beta}) is said to be {\em $d$-multiply equidistributed} in $[0,1]$ if there exists a measurable subset of seeds $T_d$, of full measure (or equivalently, $\mes(T_d) = \mes(G)$), such that the
sequence of vectors $\bbeta_k^d$  in (\ref{ex1.5mult}) is simply equidistributed
in $[0,1]^d$ for each $t\in T_d$. \\
(ii) If $\{\beta_k : k \in {\bbN}\}$ is $d$-multiply equidistributed in $[0,1]$ for all $d\geq 1$, then 
$\{\beta_k : k \in {\bbN}\}$ is called 
{\em completely equidistributed}.
\end{Definition}
We shall often abbreviate completely equidistributed sequence in $[0,1]$ by {\em c.e.s} in $[0,1]$.
Occasionally we may drop the descriptive ``in $[0,1]$''.
It is natural to also consider shifts of arbitrary fixed length.
\begin{Definition}
\label{def2}
If for some measurable subset of seeds $T_{d,h}$, again of full measure, we have that the
sequence of vectors $\bbeta_k^{d,h}$ in (\ref{ex1.6mult}) is simply equidistributed
in $[0,1]^d$ for all $t\in T_{d,h}$ then we say that  $\{\beta_k : k \in {\bbN}\}$ defined in (\ref{D:beta}) is 
 {\em $d$-multiply equidistributed in $[0,1]$ with respect to shift by $h$}.
\end{Definition}
Naturally, we omit ``with respect to shift...'' if $h=1$, and use attribute 
{\em simply} instead of $1$-multiply if $d=1$ and $h>1$.

\smallskip
\noindent
{\em Remark 1}.
(a) Fix some $h\in \bbN$ and $o\in \bbN_0$.
It is easy to check that if  $\{\beta_k : k \in {\bbN}\}$ is $d$-multiply equidistributed in $[0,1]$ with respect to shift by $h$, then the same property holds (on the same subset of seeds) for $\{\beta_{k+o} : k \in {\bbN}\}$. For this reason the definitions above are stated only in the case $o=0$.\\
(b) Due to the just mentioned easy fact,  one can quickly check (by averaging over $h$ different offsets) that if  $\{\beta_k : k \in {\bbN}\}$ is $d$-multiply equidistributed in $[0,1]$ with respect to shift by $h \geq 2$, then it is also 
$d$-multiply equidistributed in $[0,1]$ (in the sense of Definition \ref{def1}).
Note however the following example: suppose $(\alpha_k)_k$ and $(\gamma_k)_k$ are (simply) equidistributed in $[0,1]$ so that $(\alpha_k/2)_k$ and $((1 + \gamma_k)/2)_k$ are equidistributed in $[0,1/2]$ and $[1/2,1]$, respectively. Then $(\beta_k)_k$ defined by
$$
\beta_{2k-1}:=\alpha_k,\ \beta_{2k}:=\gamma_k, \quad k\geq 1,
$$
is (simply) equidistributed in $[0,1]$, but not even simply equidistributed in $[0,1]$ with respect to shift by $2$.
One could construct similar examples in the multiply equidistributed setting.\\
(c)
As already noted, the interest of including shifts becomes apparent if one compares the  equidistribution with the law of large numbers (LLN), or with Monte-Carlo simulations using pseudo-random numbers. 
If $X_1,X_2,\ldots$ is a sequence of i.i.d.~uniform random variables, and $f:[0,1]^d\to \bbR$ a bounded measurable map, then $E f(X_1,\ldots,X_d)$ is the theoretical (LLN) limit (in the $L^p$ and in the almost sure sense) of 
\begin{equation}
\label{ELLN}
\frac{1}{n} \sum_{k=1}^n f(X_{(k-1) d +1},\ldots, X_{(k-1) d +d}).
\end{equation}
Given a pseudo-random sequence $a=(a_k)_{k\geq 1}$, a direct analogue of (\ref{ELLN}) is
\begin{equation}
\label{ELLNpseudo}
\frac{1}{n} \sum_{k=1}^n f(a_{(k-1) d +1},\ldots, a_{(k-1) d +d}),
\end{equation}
and the shift by $h=d$ is typical in doing simulations.  
While it is also true that for each choice of $h,o\in \bbN$
$$
\lim_n \frac{1}{n} \sum_{k=1}^n f(X_{(k-1) h+o},\ldots, X_{(k-1) h +d +o-1}) = E f(X_1,\ldots,X_d),
$$
the variance (and therefore the theoretical error) of the approximation is the smallest if $h\geq d$. 
This variance (error) bound is constant over shifts $h\geq d$ and offsets $o\geq 1$, hence
(\ref{ELLN}) is an ``economical'' approximation (each element of $(X_k)_k$ is used once and only once), and (\ref{ELLNpseudo}) is its direct analog.\\ 
(d) We can point the reader to at least two different derivations (\cite{knu69} Section 3.5, Theorem C and the final Note in \cite{kor92} Ch.~III, \S 20) of the following important fact:  if  $\{\beta_k : k \in {\bbN}\}$ is completely equidistributed, then it is also (``omega-by-omega'', on the same set of seeds)
$d$-multiply equidistributed in $[0,1]$ with respect to shift by $h$ for each $d\geq 1$ and each $h \geq 2$. 
Due to this fact, and the observations made in (b), defining completely equidistributed with respect to shifts is superfluous. \\
(e) As a consequence of (d) and various other properties derived in \cite{F2}, Knuth \cite{knu_paper,knu69} concludes  that any c.e.s.~passes numerous empirical tests (see \cite{knu69}, Section 3.5, the comment following Definition R1), and is therefore an exemplary pseudo-random sequence.
\endpf

\smallskip
The Weyl sequence $(\beta_k)_{k\geq 1}$ from (\ref{ex1.1}) is $p$-multiply equidistributed, but it is not $p+1$-multiply equidistributed;
the multiplicatively generated sequence from (\ref{ex1.2}) is only simply equidistributed (see also Remark 2(c)).
Section \ref{S:Twonew} serves to point out that already in the linear setting, numerous examples of c.e.s.~that generalize the Weyl sequence in a natural way, can be easily constructed.

\medskip
Koksma's numbers from (\ref{ex1.3}) generally serve (see e.g.~\cite{knu69,kuni,tribble_owen,strpor}) as the prototype of a metric c.e.s.
Most of the sequel is organized in connection to this example.
In particular, Section \ref{S:Koksma} gives a page-long proof of their complete equidistribution, as an extension of the technique from Section \ref{S:linear} to the non-linear setting.
The discussion in Section \ref{S:discNT} serves to put this into perspective with respect to \cite{koksma,F2} and \cite{nietic}.
Section \ref{S:conclOP} recalls the main observations made in this and the next section, and discusses several natural open problems.

\section{Equidistribution via probabilistic reasoning}
\label{S:equi_prob}
The purpose of this section is to describe the approach of \cite{comed}, as well as to  put it into perspective with respect to \cite{koksma,F2,nietic}. 
\subsection{A lesson from the linear case}\label{S:linear}
Consider the set of multi-indices ${\bf m} = (m_1, m_2, \ldots, m_d)\in {\bbZ}^d$. 
Let $\{ x_k : k \in {\bbN}\}$ be a sequence of functions (soon taken to be linear)
as in Definitions \ref{def1}-\ref{def2},
and let points $\bbeta_k \in D$
be defined by (\ref{ex1.4mult}). 

Define the sequence $(\nu_N)_N$ of purely atomic finite (probability) measures on $[0,1]^d$ via
\begin{equation}
\label{D:nu}
\nu_N := \sum_{k=1}^N \frac{1}{N} \delta_{\bbeta_k}, \quad N\geq 1,
\end{equation}
where the dependence in $\bt\in G$ is implicit.
Then clearly, the (simple) equidistribution in $[0,1]^d$ (\ref{EeqWeyl}) is equivalent to the weak convergence of  $(\nu_N)_N$ to the uniform law on $[0,1]^d$ (denoted in (\ref{EeqWeyl}) by $\mes$), as $N$ goes to $\infty$.
This in turn is equivalent to saying that for any polynomial $p$ in $d$-variables we have $\lim_N\int p(\bs)\nu_N(d\bs)  = \int_{[0,1]^d} p(\bs)\mes(d\bs)$.

From an analyst's perspective, choosing the class of polynomials in the above characterization of weak convergence is suboptimal. Indeed, if $d=1$, the class of complex exponentials $t\mapsto \{\exp(2\pi i \,m\cdot t)\}$, $m\in\bbZ$, is orthogonal (even orthonormal) in $L^2[0,1]$, and (due to the Stone-Weierstrass theorem) the algebra they generate is
 dense in the (periodic) continuous functions on $[0,1]$, with respect to the usual sup-norm. 
The same is true in the $d$-multiple setting, this time with respect to the class of multi-dimensional complex exponentials 
$D\ni \bs\mapsto \{\exp(2\pi i \,\bm\cdot \bs)\}$, $\bm\in\bbZ^d$, where $\cdot$ is the dot (or scalar) product.
 
In order to check that $\lim_N \nu_N  = \mes$ in law, it is necessary and sufficient that for each $\bm\in \bbZ^d$
\[
\frac{1}{N} \sum_{k=1}^N \exp(2\pi i \,\bm\cdot \bbeta_k)
=
\int_D \exp(2\pi i \,\bm\cdot \bs) \nu_N(d\bs) 
\]
converges as $N\to \infty$ to
\[
\int_{[0,1]^d}\exp(2\pi i \,\bm\cdot \bs) \mes(d\bs) = \left\{
\begin{array}{ll}
1, & \bm={\bf 0},\\
0, & \bm \neq  {\bf 0}.
\end{array}
 \right.
\]
The just obtained characterization for equidistribution of $(\beta_k)_k$ is the well-known {\em Weyl criterion} \cite{We}: consider the quantities
\begin{equation}
\label{EWN} W_N(\mbox{\boldmath $\beta$},{\bf m}) \ := \ \frac{1}{N} \:\sum_{k =1}^N
 \: \exp \big ( 2 \pi i\:{\bf m} \cdot \mbox{\boldmath $\beta$}_k \big), \end{equation}
then $(\bbeta_k)_k$ is (simply) equidistributed in $D$ if and only if $W_N(\mbox{\boldmath $\beta$}, {\bf m})
\to 0$ for $N \to \infty$ and each ${\bf m} \ne {\bf 0}$.
Recalling that in our setting each $\bbeta_k$ and therefore $\nu_N$ is in fact a (measurable) function of $\bt$, the criterion reads
\begin{equation}
\label{EWeyl}
\lim_N W_N(\bbeta, \bm)= 0, \mbox{ almost everywhere in }G\mbox{, for each }\bm\in \bbZ^d\setminus\{ {\bf 0}\}.
\end{equation}

Denote by $\bbR^d\ni{\bf z} \mapsto e({\bf m},{\bf z}) = \exp (2 \pi i\, {\bf m}
\cdot {\bf z})$. 
Due to (\ref{D:beta}), and the periodicity of the complex exponential, we have the identity
\begin{equation}
\label{EidenW}
W_N(\bbeta,\bm) \equiv \frac{1}{N}\sum_{k=1}^N e(\bm,\bx_k).
\end{equation}
A crucial point is that if $({\bf x}_k)_k$ is defined by (\ref{ex1.4mult}), where the sequence of real (measurable) functions $(x_k)_k$ is given either in (\ref{ex1.1}) or in (\ref{ex1.2}), then the 
functions $G^d\ni {\bf t} \mapsto e({\bf m},{\bf x}_k({\bf t}))$ form a bounded orthogonal sequence in $L_2(G,\bbC)$.
Note that this is the first time that the linearity in $t$ (resp.~$\bt$) of the functions in (\ref{ex1.1},\ref{ex1.2})
(resp.~(\ref{ex1.4mult})) is being called for.

Unlike \cite{comed}, we continue the discussion
using probabilistic wording and notation.
Given $z\in \bbC$, denote by $\overline{z}$ its complex conjugate. 
For a fixed $\bm\in  \bbZ^d \setminus \{0\}$ and each $k\in \bbN$, define $Y_k(\bm\,;\bt):= e(\bm,\bx_k(\bt))$,  $\bt\in G$. 
Then  $(Y_k(\bm))_k$ is a sequence of complex-valued random variables on the probability space $(G,{\mathcal B}, P)$, where ${\mathcal B}$ is the Borel $\sigma$-field on $G$, and $P(d\bt)=\mes(d\bt)/\mes(G)$, such that $E Y_k(\bm)=0$, for all $k\geq 1$, and
\begin{equation}
\label{Euncorrel}
E(Y_k(\bm) \overline{Y_l(\bm)})=  \left\{
\begin{array}{ll}
1, & k=l,\\
0, & k\neq l,
\end{array} 
\right.\quad l,k\in \bbN.
\end{equation}
The standard proof of the strong LLN for i.i.d.~variables with finite second moment (see for example the exercise concluding \cite{durrett} Ch. I, Section 7) can clearly be extended to sequences of pairwise uncorrelated centered random variables with constant (or uniformly bounded) variance.
Anticipating some readers with non-probabilist background, we include a sketch of the argument in points LLN.(a)-LLN.(d) below.
 The above sequence $(Y_k(\bm))_k$ has the just stated properties, and therefore if $S_N\equiv S_N(\bm):=\sum_{k=1}^N Y_k(\bm)$, then
\begin{equation}
\label{Econcl}
\lim_N \frac{S_N}{N} = \lim_N \frac{1}{N} \sum_{k=1}^N Y_k(\bm) = 0,\quad \mbox{ almost surely}.
\end{equation}
Recalling (\ref{EidenW}) and the fact that $\bm\in \bbZ^d \setminus \{0\}$ was fixed but arbitrary, gives (\ref{EWeyl}) ($\mes$-a.e.~is the same as $P$-a.s,) and therefore the above stated (simple) equidistribution in $[0,1]^d$
of the sequence of vectors $(\bbeta_k)_k$.

\smallskip
\noindent
LLN.(a) 
Note that $ES_N=0$ and ${\rm var}(S_N)=O(N)$.\\
LLN.(b) 
Use the Chebyshev (or the Markov) inequality, and the Borel-Cantelli lemma on the subsequence $N_n= n^2$ to conclude that 
$$ \sum_{n=1}^\infty P( |S_{n^2}| > n^2 \eps) = \sum_{n=1}^\infty  O(1/n^2)<\infty, \quad \forall \eps>0,$$
and therefore that $S_{N_n}/{N_n}\to 0$, almost surely.\\
LLN.(c) 
Since 
$E(\sum_{N\in   [n^2+1,(n+1)^2]} |S_N - S_{n^2}|^2) = \sum_{N\in   [n^2+1,(n+1)^2]}  O(N-n^2) = O(n^2)$,
and therefore $E(\max_{N\in   [n^2+1,(n+1)^2]} |S_N- S_{n^2}|^2) \leq D n^2$, another application of the Markov inequality gives that with probability greater than $1- \frac{2D}{\eps^2 n^2}$
\begin{equation}
\label{Esandwich} 
|S_N - S_{n^2}| \in [ - \eps n^2 , \eps n^2], \quad \forall N\in[n^2+1,(n+1)^2-1],
\end{equation}
and again by the Borel-Cantelli lemma, that (\ref{Esandwich}) happens with probability $1$ for all but finitely many $n$.\\
LLN.(d) Divide (\ref{Esandwich}) by $n^2$, noting that $N=n^2 + O(n)$.
Use the fact $\eps>0$ was arbitrary,
and the conclusion of (b).

\smallskip
\noindent
{\em Remark 2}.
(a) In our special setting, each member of the sequence $(Y_k)_k$ is uniformly bounded below (and above), so (\ref{Esandwich}) could have been replaced by a simpler estimate: with probability $1$ for all $n$
$$
|S_N - S_{n^2}| \in [- (N-n^2)c,  ((n+1)^2-N)c], \quad N\in[n^2+1,(n+1)^2-1].
$$
(b) If $r=d=1$, then (\ref{ex1.4mult}) and (\ref{ex1.5mult}) coincide, leading to the conclusion that if
 $(x_k)_k$ is again either from (\ref{ex1.1}) or in (\ref{ex1.2}), then the corresponding $(\beta_k)_k$ is simply equidistributed.
The study of multiple equidistribution can be similarly set up (see also Section \ref{S:Twonew} below): here (\ref{D:nu}) has the same form, but we take $d\geq 2$ and $\bbeta:=(\bbeta_k^d)_k$ defined as in (\ref{ex1.5mult}) with $(x_k)_k$ linear (as in (\ref{ex1.1},\ref{ex1.2})), and study the averaged Weyl sums of terms
$Y_k(\bm\,;t)=e(\bx_k(t),\bm)$, for $t\in G$, for each $\bm \in \bbZ^d\setminus \{\bf 0\}$.\\
(c)
Recall the multiplicatively generated sequence from (\ref{ex1.2}). To see that it
is not $2$-multiply equidistributed, take $\bm=(M,-1) \neq (0,0)$ and note that with this choice of $\bm$ the corresponding criterion (\ref{EWeyl}--\ref{EidenW}) converges nowhere to $0$. Indeed, $m_1 x_k + m_2 x_{k+1} \equiv 0$, for any $k\in \bbN$.\\
Similarly, it is possible to find a non-trivial $p+1$-dimensional vector $\bm$ such that 
\[
m_1 k^p  + m_2 (k+1)^p + \ldots + m_p (k+p-1)^p + m_{p+1} (k+p)^p \mbox{ is constant over }k. 
\]
This amounts to solving a linear system $A \bx = {\bf 0}$, where $A$ has $p$ rows and $p+1$ columns.
The corresponding criterion (\ref{EWeyl}--\ref{EidenW}) applied to Weyl numbers (\ref{ex1.1}) will again converge nowhere to $0$ (non-convergence to $0$ on an event of positive probability would already be sufficient).

\subsubsection{Equivalent formulations of complete equidistribution}
\label{S:equiv_form}
Let $d\geq 1$ be fixed.
The $d$-dimensional (or multiple) {\em discrepancy}  at level $N$ is defined by
$$ D_N^d(\bt):= \sup_{A\in J}
  \left|  \frac{\sum_{k = 1}^N \:{\bbJ_A(\bbeta_k^d(\bt))}}{N} - \mes(A)  \right|, \quad \bt\in G,
$$
where $J$ is the family of boxes as in definition (\ref{EeqWeyl}).
The $d$-dimensional {\em ``star'' discrepancy}  at level $N$ is defined by
$$ D_N^{*,d} (\bt):= \sup_{A\in J^*}
  \left|  \frac{\sum_{k = 1}^N \:{\bbJ_A(\bbeta_k^d(\bt))}}{N} - \mes(A)  \right|, \quad \bt\in G,
$$
where $J^*$ is the family of boxes in $[0,1]^d$ as above, with lower left corner equal to ${\bf 0}$. 
Equivalently, the supremum above is taken over all boxes $A$ of the form $\prod_{i=1}^d (0,b_i]$, where $0\leq b_i\leq 1$, for each $i$.
Note that we again deviate slightly from the standard definitions (see Kuipers and Niederreiter \cite{kuni}), where discrepancy sequences are defined for deterministic sets or sequences. Easy properties of measurability of functions make each $D_N^d$ and $D_N^{d,*}$ a random variable in the current (stochastic a.s./metric) setting. Note that $D_N^{*,d}$ is a direct extension of the {\em Kolmogorov-Smirnov} statistic to the equidistribution setting.

Recall (\ref{D:nu}) with $\bbeta_k=\bbeta_k^d$, $k\geq 1$.
As any probabilist knows (think about {\em cumulative distribution functions}), the weak convergence  of the sequence of (random) measures $\nu_N$
to the uniform law on $[0,1]^d$ is, omega-by-omega (where ``omega'' is typically denoted by $\bt$), equivalent to the convergence of  $(D_N^d)_N$ (or equivalently of $(D_N^{*,d})_N$) to $0$ (a clear generalization of the Glivenko-Cantelli theorem). 
Let us paraphrase: $(\beta_k)_k$ is $d$-multiply equidistributed in $[0,1]$ if and only if $\lim_N D_N^d=\lim_N D_N^{*,d}=0$, almost surely.
In particular, a sequence $(\beta_k)_k$ is a c.e.s.~in $[0,1]$ if and only if 
\begin{equation}
\label{Eequiv_form_one}
\lim_N D_N^d=\lim_N D_N^{*,d}=0, \ \forall d\in \bbN, \ \mbox{almost surely}.
\end{equation}

Note that $I_{\infty}:=[0,1]^\bbN$ is a product of compact spaces, and therefore compact itself.
Let ${\mathcal B}_m$ be the Borel $\sigma$-field on $[0,1]^m$.
Consider the {\em cylinder sets} 
$C\subset I_{\infty}$, 
such that $C\=C_b\times [0,1]^{\bbN}$, for some $C_b\in {\mathcal B}_m$ and $m\in \bbN$.
Let ${\mathcal B}_{\infty}$ be the $\sigma$-field on $I_\infty$ generated by the cylinders.
Instead of $d$-dimensional vectors, one could consider straight away the ``$\infty$-dimensional'' (random) vector sequence
$$
\bx_k(\bt)=(x_k(\bt),x_{k+1}(\bt),x_{k+2}(\bt),\ldots), \ k\geq 1,\ \bt \in G,
$$
and its corresponding $\bbeta_k^\infty:=  \bx_k \modb$, where again ``modulo'' operation is applied component-wise (a.s.).
Let $\nu_N^\infty$ be as $\nu_N$ in (\ref{D:nu}), but with $\bbeta$ redefined as $\bbeta^\infty$.
Another formulation of complete equidistribution can be read off from an ``abstract fundamental theorem'' (see e.g.~\cite{kuni}, Ch.~3, Theorem 1.2 and the remark following it) or easily proved by approximating all closed sets in ${\mathcal B_\infty}$ by closed cylinders: $(\beta_k)_k$ is a c.e.s.~in $[0,1]$ if and only if $(\nu_N^\infty)_N$ converges weakly to the uniform law  on $I_\infty$. A realization $\Gamma$ from this limiting uniform law (that is, the random object having that law) is also called the {\em infinite statistical sample} or the  {\em i.i.d.~family of uniform $[0,1]$ random variables} $(U_1,U_2,\ldots)$.
The just made statement could also be reformulated as follows:
$(\beta_k)_k$ is a c.e.s.~in $[0,1]$ if and only if 
for any $f$ bounded and continuous function on $I_\iy$ (equipped with the product topology), we have
\begin{equation}
\label{Eequiv_form_Gamma}
\lim_N\frac{1}{N}\sum_{k=1}^N f(\bbeta_k^\infty) = {\bf{E}}\, (f(\Gamma)), \ \mbox{almost surely}.
\end{equation}

As indicated in Section \ref{S:lite}, another formulation/criterion of complete equidistribution for a (deterministic) sequence of real numbers in $[0,1]$ was derived by Knuth \cite{knu_paper,knu69}. There is no doubt that this could  also be turned into a stochastic (a.s.) formulation for a sequence (\ref{D:beta}), the details are left to an interested reader. 

\subsubsection{Why not simply take an i.i.d.~family of uniforms?}
\label{S:why_not}
The first author can guess that, especially on the first reading, a non-negligible fraction of fellow probabilists 
could be asking the above or a similar question. 
It is clear that an i.i.d.~sequence $\Gamma$ of uniform random variables is a c.e.s.~in $[0,1]$ in the sense of Definition \ref{def1}, or any of its equivalent formulations described in the previous section. 
Studying complete equidistribution, starting from i.i.d.~(or similar) random families is precisely what probability oriented works \cite{hol1,hol2,lacaze,loynes} did.
However, to most non-probabilist mathematicians this will not mean much, especially due to the fact that the rigorous probability theory is axiomatic, and the rigorous construction of $\Gamma$ rather abstract.

More precisely, $\Gamma$ is not presented in a neat (classical) functional format, like any of the sequences $(x_k)_k$ and their corresponding $(\beta_k) = (x_k \mod)_k$ in (\ref{ex1.1},\ref{ex1.2},\ref{ex1.3}). Instead we usually start with an abstract infinite product space  (or more concretely, with $[0,1]^\bbN$), and the $k$th uniform random variable equals the identity map from the $k$th space to itself. This suffices for the most purposes of modern probability theory, but may not seem very convincing to a non-probabilist who wishes to ``see a concrete example of'' $\Gamma$. A related wish to ``see an explicit outcome from a concrete example of $\Gamma$'' quickly leads to philosophical discussions around the question ``What is a random sequence?'' (the reader is referred to the section in \cite{knu69} carrying that very title). 

The Kac \cite{kac} approach to the i.i.d.~discrete-valued sequence is also revolutionary from the point of view of the just mentioned drawback.
Indeed, an infinite sequence of independent Bernoulli$(1/2)$ random variables (also called the {\em Bernoulli scheme}) can be explicitly constructed on $\Omega=[0,1]$, with $\FF$ equal to the Borel $\sigma$-field ${\mathcal B}$, and $P$ equal to $\lambda$.
This was done in \cite{kac}, by applying a simple transformation to the classical Rademacher system. We leave on purpose the link with (\ref{ex1.2}), and other details, to interested readers, as well as the discovery of the related $b$-adic Rademacher system and its relation to the discrete uniform law on $\{0,1,\ldots,b-1\}$.

Given an infinite sequence $X:=(X_i)_{i=1}^\iy$, which has the Bernoulli scheme distribution (take for example the one from the above recalled construction by Kac), one can define $\Gamma$ on the same $(\Omega,\FF,P)$ by ``redistributing the digits'' via a triangular scheme, for example $U_1:= X_1/2 + X_3/2^2 + X_6/2^3+X_{10}/2^4+\cdots$, $U_2:= X_2/2+X_5/2^2+X_9/2^3+\cdots$, $U_3:=X_4/2+X_8/2^2+\cdots$, $U_4:=X_7/2+\cdots$,  and so on, and finally $\Gamma=(U_i)_i$.
Is this asymptotic definition of $\Gamma$ on $([0,1],{\mathcal B}, P)$ sufficiently explicit?
Or are the examples like (\ref{ex1.3},\ref{ex_novi1},\ref{ex_novi2}) - all leading (after $\!\!\mod$ application) to c.e.s.~but none to i.i.d.~uniforms - more reassuring to think about (more amenable to analysis)? The answer will likely vary from one peer to another, depending not only on their mathematical background and research interests, but also on their personal perception of randomness. Note that such and related questions challenged the very founders of probability theory only 80 years ago (see e.g.~the historical notes of Knuth \cite{knu69}, Section 3.3.5, or \cite{vanLambalgen,burdzy,burdzy_res}).

For reader's benefit,  we mention in this paragraph two other prominent approaches to  randomness, nowadays practically forgotten by the probability and statistics communities: the von Mises collectives, and Martin-L\"of sequences.
The original definition \cite{vonMises}  of {\em Kollektiv} (or {\em collective}) is not transparent. Following \cite{burdzy_res} (e.g. page 49),
a collective is an infinite sequence of observations such that the relative frequency of an event converges to the same number along every subsequence chosen without prophetic powers; this common limit is called the probability of the event in the collective.
For further interpretations, discussions and more see \cite{vanLambalgen,burdzy,burdzy_res}.
Continued search for (more concrete) examples of mathematically generated randomness led to {\em Kolmogorov's complexity} \cite{kolmogorov2}, and the {\em Martin-L\"of sequences} (see \cite{MartinLof} and \cite{levinL}).
Informally speaking,
Martin-L\"of \cite{MartinLof} proves existence of binary sequences that pass the  "super-test for randomness".
Any such sequence can be rightfully called random.  Concreteness of this definition is again relative, and therefore left to an individual appraisal (Burdzy gives an interesting commentary in \cite{burdzy_res}, pp.~77-78).

\subsubsection{Novel classes of examples of c.e.s.}
\label{S:Twonew}
Define 
\begin{equation}
\label{ex_novi1}
 x_k(t) \ = \ k^k\,t, \quad t \in G \::=\:(0,1),\ k\geq 1.
\end{equation}
and
\begin{equation}
\label{ex_novi2}
 x_k(t) \ = \ k!\,t, \quad t \in G \::=\:(0,1),\ k\geq 1.
\end{equation}
We claim that (\ref{D:beta}) defined using either (\ref{ex_novi1}) or (\ref{ex_novi2})
is completely equidistributed in the sense of Definition \ref{def1}.  

Let us take (\ref{ex_novi2}) for example.
Simple equidistribution can be quickly deduced as in the previous section (see Remark 2(b) in particular).
Now fix some $d\geq 2$ and $\bm\in \bbZ\setminus \{\bf 0\}$. 
WLOG we can assume that the Weyl criterion (\ref{EWeyl}) has been shown for the $(d-1)$-multiply case.
Moreover, due to symmetry, we can and will assume $m_d>0$.
We have
\begin{equation}
\label{D:Yk}
Y_k(t)\equiv Y_k(\bm\,;t) := e(\bm, \bx_k(t)), 
\end{equation}
which in this special case reads 
$Y_k(t)=\exp(2 \pi i  \,t\sum_{i=1}^d m_i (k+i-1)!)$, $t\in [0,1]$. 
It is clearly true that $E(Y_k)=0$, and $E|Y_k|^2 =1$, for each $k\in \bbN$.
It is furthermore easy to see that whenever $k-l > c(\bm)$ for some $c(\bm)\in \bbN$, then $\sum_i m_i [(k+i-1)! - (l+i-1)!])$ is a strictly positive integer, and due to the periodicity of sine and cosine, we have again $E(Y_k \overline{Y_l})= \overline{E(\overline{Y_k}Y_l)} =0$ (or equivalently, $Y_k$ and $Y_l$ are uncorrelated). 
Moreover, $|Y|$s are uniformly bounded by $1$.
Therefore 
$$
{\rm var}\left(\sum_{k=1}^N Y_k\right) = N + \sum_{k=1}^N \sum_{l=k-c(\bm)}^{k-1} 1 = N(1+ c(\bm)),
$$
and again the argument LLN.(b)-LLN.(d) (using the simplification from Remark 2(a)) applies.
By induction on $d$, we obtain the above stated complete equidistribution.
The proof for (\ref{ex_novi1}) follows the very same steps.
The reader is undoubtedly able to construct further new examples of c.e.s.~for which the elements of $(Y_k)_k$ are pairwise uncorrelated (except for nearby indices).

\subsection{The Weyl variant of the SLLN and generalizations}
\label{S:weylSLLN}
It is easy to see that the steps in LLN.(a)-LLN.(d) could be applied to show simple equidistribution in $[0,1]$ of (\ref{D:beta}), whenever $x_k(t):=a_k t$, $k\geq 1$ and $(a_k)_k$ is a sequence of distinct integers, yielding an alternative derivation of
\cite{kuni} I, Theorem 4.1.
This result is already due to Weyl \cite{We}, and  its original proof made a profound impact on the development of (analytic) number theory.
Indeed, in the above and analogous situations (see e.g.~\cite{kuni,MF,drmota_tichy}, as well as \cite{hol1}, Theorem 2), it is standard to apply the following variant of the SLLN, due to Weyl \cite{We}: since $E|S_N|^2 =O(N)$, then
$$
\sum_n E\left(\frac{|S_{n^2}|}{n^2}\right)^2 <\infty,
$$
and therefore by Tonelli's theorem $\sum_n \frac{|S_{n^2}|}{n^2}<\infty$ almost surely, and in particular $S_{n^2}/n^2\to 0$, as $n\to \infty$. Finally use Remark 2(a) to get (\ref{Econcl}).

The Weyl variant of the SLLN has the following important extension, 
due to Davenport, Erd\"os and LeVeque \cite{delv63}, that has been since used extensively in number theoretic studies (see e.g.~\cite{kuni} I, Theorem 4.2):
\begin{equation}
\label{Edelv_SLLN}
\mbox{if }\ \ \sum_n \frac{1}{n} E\left(\frac{|S_n|}{n}\right)^2 <\infty  \ \mbox{ then } (\ref{Econcl}).
\end{equation}
Its proof strongly relies on the uniform boundedness of the individual random variables $Y$, but not at all on their particular (complex exponential) form (see Lyons \cite{lyo88}, Theorem 1). 
Lyons \cite{lyo88} takes moreover general complex-valued sequence $(Y_k)_k$, and derives different generalizations of the above SLLN criterion of \cite{delv63}, where the uniform boundedness condition is replaced by various bounded moment conditions.

\smallskip
\noindent
{\em Remark 3}. In view of the SLLN stated here, the fact that the sequences of functions  from Section \ref{S:Twonew} generate c.e.s.~is again trivial. Just like with the usual SLLN (of LLN.(a)-(d) and Remark 2(a)), it is important here that all the covariances can be computed or adequately (uniformly) estimated. 

\subsection{Weyl-like criterion for weakly c.e.s.}
\label{S:WCUD}
Recall the (random) discrepancy sequences introduced in Section \ref{S:equiv_form}.
Owen and Tribble (e.g.~\cite{owen_tribble} Definition 2 or \cite{tribble_owen} Definition 5) define {\em weak complete equidistribution} (aka WCUD) in terms of the weak(er) convergence of the discrepancy sequence as follows: $(\beta_k)_k$ is WCUD if and only if 
$D_N^{*,d} \Rightarrow 0$, or equivalently (since the limit is deterministic), if $D_N^{*,d} \to 0$ in probability, for each $d\geq 1$.

The following ``a.s.-convergence along subsequences'' characterization is well-known: a sequence of random variables $(X_n)_n$ converges in probability to $X$, if and only if for any subsequence $(n_k)_k$, one can find a further subsequence $(n_{k(j)})_j$ such that $X_{n_{k(j)}}$ converges to $X$ almost surely.
Let us apply it to the discrepancy sequences, and conclude (due to the countability of $\bbN$, and the Arzel\`a-Ascoli diagonalization scheme) that for any  subsequence $(n_k)_k$ one can find a further subsequence $(n_{k(j)})_j$
and a single event (Borel measurable set) $G_f$ of full probability such that 
$$
D_{n_{k(j)}}^{*,d}(\bt) \to 0, \ \mbox{ for all } d \geq 1, \mbox{ and all } \bt \in G_f \subset G. 
$$
But now we recall (as in Section \ref{S:equiv_form}) that the above convergence is equivalent ($\bt$-by-$\bt$) to the simultaneous (for each $d\in \bbN$ and on the full probability event $G_f$) weak convergence of the random sequence $(\nu_{n_{k(j)}}^d)_j$ 
to the corresponding uniform law $\mes$ on $[0,1]^d$, where
$$
\nu_{N}^d := \sum_{k=1}^N \frac{1}{N} \delta_{\bbeta_k^d},  \ d\in \bbN.
$$

The last made claim is in turn equivalent (recall the reasoning of Section \ref{S:linear} and definition (\ref{D:Yk})) to the statement
\begin{equation}
\label{EWeyl_subseq}
\lim_j \frac{1}{n_{k(j)}}\sum_{k=1}^{n_{k(j)}}Y_k(\bm)= 0, \mbox{ on }G_f\mbox{, for each }\bm\in \bbZ^d\setminus\{ {\bf 0}\},
\mbox{ and each } d\in \bbN.
\end{equation}
This leads to the following conclusion: for each subsequence $(n_k)_k$ one can find a further subsequence $(n_{k(j)})_j$ such that for each $\bm\in \bbZ^d\setminus\{ {\bf 0}\}$ the rescaled Weyl sums indexed by $\bm$ converge to $0$  along that sub(sub)sequence, almost surely. In other words, for each $d$ and $\bm \in \bbZ^d\setminus\{ {\bf 0}\}$, 
instead of (\ref{Econcl}) we arrived to 
\begin{equation}
\label{Econcl_proba}
\lim_N \frac{S_N}{N} = \lim_N \frac{1}{N} \sum_{k=1}^N Y_k(\bm) = 0,\quad \mbox{in probability}.
\end{equation}
Finally recalling that, in our special setting, $|S_N|/N$ is uniformly bounded by $1$, the (adaptation of) Lebesgue dominated convergence theorem says that (\ref{Econcl_proba}) is equivalent to 
\begin{equation}
\label{Econcl_Lone}
\lim_N E\frac{S_N}{N} = \lim_N E\frac{|S_N|}{N}   = 0.
\end{equation}
We record the above sequentially made  equivalences as 
\begin{Lemma}[Sufficiency and necessity for weak c.e.s.~(aka WCUD)]
\label{L:crit_Wces}
$ $\\ 
Let $(x_k)_k$ be a sequence of  random variables (or measurable functions on $G$), $Y_k(\bm)$ be defined  by (\ref{D:Yk}), and  $(\beta_k)_k$ by (\ref{D:beta}).
Then $(\beta_k)_k$ is weakly completely equidistributed in $[0,1]$ if and only if 
(\ref{Econcl_Lone}) is valid for each $d$ and $\bm \in  \bbZ^d\setminus\{ {\bf 0}\}$.
\end{Lemma}
In view of Lemma \ref{L:crit_Wces}, the claim from \cite{tribble_owen}  that WCUD sequences are hard to construct seems rather surprising.

\subsection{Extensions to the non-linear setting}
\label{S:Koksma}
As already noted, the reasoning of Section \ref{S:linear} is crucially connected to linearity only through (\ref{Euncorrel}). In particular, if we were to take (\ref{ex1.3}) or some other sequence of functions, then apply (\ref{ex1.5mult}) to obtain $d$-dimensional vector sequences, and plug them into (\ref{EWN},\ref{EidenW},\ref{D:Yk}), the criterion (\ref{EWeyl}), to be verified through (\ref{Econcl}), would stay the same. 
Note however that the correlation structure may (and typically does) become much more complicated than that given in (\ref{Euncorrel}) or in Section \ref{S:Twonew}.
For a detailed discussion of problems with non-linearity see \cite{aistleitner_MG}.

In the special case of Koksma's numbers, the probability $P$ is set to $\mes$ renormalized by $(a-1)$ on $[1,a]$, or equivalently, $P$ is the uniform law on $[1,a]$. 
Before considering covariance, one could already note that $Y_k$ may not (necessarily) be even centered.
Since $Y_k(t)\equiv Y_k(\bm\,;t):= e(\bm,\bx_k(t))$, $t\in [1,a]$, we have that
\begin{equation}
\label{EexpYk}
E(Y_k) =\frac{1}{a-1} \left[\int_{1}^a \cos(2\pi \bm \cdot \bx_k(t))  \,dt + i  \int_{1}^a \sin(2\pi \bm \cdot \bx_k(t))  \,dt\right].
\end{equation}
The non-zero expectation does not matter much, if one could show that the long term average of $(Y_\cdot)$ is approximately centered, or equivalently, that
$$
\lim_N \frac{1}{N}\sum_{k=1}^N EY_k=0.
$$
The covariance analysis to be done below will imply an even stronger estimate $E(|S_N|)= O(N^{1-1/d})$, but already from Jensen's inequality we get that if $E(|S_N|^2) = o(N^2)$, then $E|S_N|=o(N)$, which would be sufficient.

The argument LLN.(a)-LLN.(d) will continue to work for the sequence $(Y_k-EY_k)_k$, and  
will imply (\ref{Econcl}) for $(Y_k)_k$, even in the presence of correlation, provided that the total covariance    
$$
\sum_{k=1}^n \sum_{l=k+1}^n |E(Y_k \overline{Y_l})+ E(Y_l \overline{Y_k})|
$$ 
grows sufficiently slowly in $N$.
Indeed, a little thought is needed to see that if
\begin{equation}
\label{Einegkl}
|E(Y_k \overline{Y_l}+ Y_l \overline{Y_k})| = \frac{1}{a-1}\left|\int_{[1,a]} \cos(2\pi \, \bm \cdot (\bx_k(t)-\bx_l(t)) \,dt\right| \leq c|k-l|^{-\delta}, \ k\neq l,
\end{equation}
for some $c=c(a,\bm,d)<\infty$, $\delta=\delta(a,\bm,d)>0$, then 
the conclusion (\ref{Econcl}) remains (this time choose $N_n=\lfloor n^{2/\delta}\rfloor$, and then apply the sandwiching argument).
Alternatively, this can be immediately deduced form the SLLN criterion (\ref{Edelv_SLLN}).

We record the just made observations in form of a lemma.
\begin{Lemma}[Sufficiency for complete equidistribution I]
\label{L:crit_covarianceI}
Let $(x_k)_k$ be a sequence of  random variables (or measurable functions on $G$), $Y_k(\bm)$ be defined  by (\ref{D:Yk}), and  $(\beta_k)_k$ by (\ref{D:beta}).
If for each $\bm$,
 $|E(Y_k(\bm) \overline{Y_l}(\bm)+ Y_l(\bm) \overline{Y_k}(\bm))| = O(|k-l|^{-\delta})$, for some $\delta(\bm)>0$, then $(\beta_k)_k$ 
 is completely equidistributed in $[0,1]$.
\end{Lemma}

Deriving (\ref{Einegkl}) is a longer calculus exercise, sketched here for reader's benefit.
Fix $d\geq 1$, $\bm\in \bbZ^d \setminus \{\bf 0 \}$, and consider $\bbeta_k^{d}= (\beta_{k1},\ldots,\beta_{kd})$ defined in (\ref{ex1.5mult}).
Recalling that $x_k(t)= t^k$ we have
 $\bm \cdot (\bx_k(t)-\bx_l(t)) =
 \sum_i m_i (t^{k+i-1} -  t^{l+i-1})$. If $p(t)\equiv p(t;\bm)= \sum_i m_i t^{i-1}$ and  $g(t)\equiv g(t;k,l)=t^{k} - t^{l}$, then
\begin{equation}
\label{Epg}
\bm \cdot (\bx_k(t)-\bx_l(t))  = p(t) g(t), \quad t\in G=[1,a].
\end{equation}

One can assume WLOG that $k>l$, so that $g(t)$ is non-negative on $G$, with a single zero $t=1$.
Let $d_*$ equal to the maximal index $i$ such that $m_i\neq 0$.
The polynomial $p$ does not depend on $k,l$. It has degree $d_*-1$ and therefore at most $d_*-1\leq d-1$ real zeros, each of which may fall into $G$.
Therefore $p(\cdot) g(\cdot)$ takes value $0$ at $1$, and at most $d-1$ other points in $G$. Let us denote these zeros by $(z_j)_{j=1}^{r}$, where $1=:z_0 \leq z_1<\ldots<z_r\leq a=:z_{r+1}$. 
For $j=1,\ldots, r$, let $l_j$ be the multiplicity of $z_j$ for $p$. Extend this definition to  $l_0=0$ if $z_1>1$, and  $l_{r+1}=0$ if $z_r<a$.
Define
$$
b_1:=\max_{t\in [1,a]} |p'(t)|,\quad b_2:=\max_{t\in [1,a]} |p''(t)|.
$$

Keeping in mind that $r\leq d-1$,
we now let $I_i:=[z_i,z_{i+1}]$, and  estimate
\begin{equation}
\label{Eint_i}
\left|\int_{I_i} \cos(2\pi \, g(t)p(t))\,dt\right|
\end{equation}
separately for each $i=0,\ldots,r$. 
The non-zero polynomial $p$ does not change sign on $I_i$, so WLOG we can assume that it takes positive values in the interior of $I_i$.
Moreover 
\begin{equation}
\label{Ebdp}
\min_{z\in I_i} p(z)/((z-z_i)^{l_i} (z_{i+1}-z)^{l_{i+1}})=:c_i>0.
\end{equation}
If $z_{i+1}-z_i \leq\frac{3}{(k-l)^{1/d}}$ then (\ref{Eint_i}) is clearly bounded by $3/(k-l)^{1/d}$.
Otherwise, define $I_i':= [z_i + \frac{1}{(k-l)^{1/d}}, z_{i+1} - \frac{1}{(k-l)^{1/d}}]$.
Split the domain of integration in (\ref{Eint_i}) into three disjoint pieces:
$[z_i, z_i + \frac{1}{(k-l)^{1/d}})$, $I_i'$ and $(z_{i+1} - \frac{1}{(k-l)^{1/d}}, z_{i+1}]$.
The integral of $\cos(\cdot)$ over the first and the final piece is again trivially bounded above by $1/(k-l)^{1/d}$. 
For the middle piece, it suffices to show that:

\vspace{0.1cm}
(a) $(pg)'(z_i + \frac{1}{(k-l)^{1/d}})\geq \bar{c}_i |k-l|^{1/d}$, where $\bar{c}_i>0$ and that

\vspace{0.05cm}
(b) $(pg)''(t) = g''(t) p(t) + 2 g'(t)p'(t) + g(t) p''(t)>0$  on $I_i'$.

\vspace{0.1cm}
\noindent
Indeed, (a)-(b) would imply that $2 \pi g(\cdot)p(\cdot)$ increases more and more rapidly on $I_i'$, and in turn that $ \cos(2\pi \, g(t)p(t))$ makes shorter and shorter ``excursions'' of alternating sign, away from $0$.
This would yield an upper bound for the integral over $I_i'$ in terms of  $t_i^*- (z_i + \frac{1}{(k-l)^{1/d}})$ where $t_i^*$ is the second zero  of  $ \cos(2\pi \, g(\cdot)p(\cdot))$ in $[z_i + \frac{1}{(k-l)^{1/d}}, \infty]$.
From (a) and (b) one could easily conclude that $t_i^*- (z_i + \frac{1}{(k-l)^{1/d}})\leq 
(2\bar{c}_i)^{-1}(k-l)^{-1/d}$.
The just made reasoning justifies the upper bound 
$(2+(2\bar{c}_i)^{-1}) (k-l)^{-1/d}$, which 
 in turn implies the upper bound (\ref{Einegkl}) with $\delta=1/d$ and $c(a,\bm,d):=\sum_{i=0}^r(2+(2\bar{c}_i)^{-1})$.
Alternatively, one could use \cite{kuni} I, Lemma 2.1.

To show (b), note that $g''(t) = k(k-1) t^{k-2} - l(l-1)  t^{l-2}$ is greater 
than  
$(k^2 - l^2 + O(k+l)) t^{k-2}$.
We also have  $g'(t) = k  t^{k-1} -  lt^{l}\leq ka \cdot t^{k-2}$ and $g(t)\leq a^2 \cdot t^{k-2}$. 
Due to (\ref{Ebdp}) and the fact $l_i+l_{i+1}\leq d_*-1\leq d-1$,  the leading term of $(pg)''(t)$ is positive and bounded below by $c_i |k-l|^{1/d}(k+l)t^{k-2}$ on $I_i'$, while the two other terms are bounded above by $b_1 k a t^{k-2}$ and $b_2 a^2t^{k-2}$, respectively, yielding (b).
We leave to the reader a similar (and easier) argument for (a).

Finally note that the derivation of the bound in (\ref{Einegkl}) applies also to estimating both the real and the imaginary part of 
(\ref{EexpYk}), and leads to an analogous bound $EY_k=O(k^{-1/d})$.

\subsubsection{Koksma's numbers have even stronger properties}
\label{S:discNT}
Koksma's derivation \cite{koksma} of the simple equidistribution of Koksma's numbers was similar to that included above (see also \cite{kuni} I, Theorem 4.3), and easier due to the fact that a simpler polynomial $m\cdot g$ replaces $p\cdot g$ in (\ref{Epg}).
As already mentioned in the introduction, Franklin \cite{F2} found a way
(again similar to the calculus exercise above)
of building on Koksma's analysis, without directly applying Weyl's criterion for multiple equidistribution. 

Koksma \cite{koksma} showed in addition a stronger type of equidistribution:
suppose that $(a_k)_k$ is a sequence of positive distinct integers, and consider the reordering (with possible deletion) of Koksma's numbers obtained via (\ref{D:beta}) from $(x_{a_k})_{k\geq 1}$ defined as in (\ref{ex1.3}), the obtained sequence is again simply equidistributed.  
Niederreiter and Tichy \cite{nietic} proved that the above arbitrarily permuted (sub)sequence of Koksma numbers is in fact completely equidistributed, thus confirming a guess of Knuth's \cite{knu69}, who called the above property pseudo-random of R4 type. 
As it turns out, the arguments in \cite{nietic} are equally based on probabilistic reasoning (using the \cite{delv63} variant (\ref{Edelv_SLLN}) of the SLLN and covariance calculations). 
For the benefit of the reader we  sketch it next. 

Suppose initially that, for each $\bm\neq {\bf 0}$,
$|E(Y_k(\bm) \overline{Y_l}(\bm)+ Y_l(\bm) \overline{Y_k}(\bm))|$ is bounded by $c(\bm)/\log^2{N}$,
unless $|k-l|\leq \sqrt{N}$ or $\min\{k,l\}\leq \sqrt{N}$.
Then it is easy to see that
$E|S_N|^2 \leq N + 2\sum_{k=1}^N \sqrt{N} + \sum_{k=\sqrt{N}}^N \sum_{l=k+\sqrt{N}}^N c(\bm)/\log^2{N} = N + N^{3/2} + N^2/\log^2{N}$ which is $O(N^2/\log^2{N})$.
The criterion (\ref{Edelv_SLLN}) implies the required convergence in the Weyl criterion.
We again record the just made observations in form of a lemma.
\begin{Lemma}[Sufficiency for complete equidistribution II]
\label{L:crit_covarianceII}
let $(x_k)_k$ be a sequence of  random variables (or measurable functions on $G$), $Y_k(\bm)$ be defined  by (\ref{D:Yk}), and  $(\beta_k)_k$ by (\ref{D:beta}).
If for each $\bm$, 
$|E(Y_k(\bm) \overline{Y_l}(\bm)+ Y_l(\bm) \overline{Y_k}(\bm))| = O(1/\log{N}^2)$, whenever $\min\{l,k-l\}\geq \sqrt{N}$, then $(\beta_k)_k$ 
 is completely equidistributed in $[0,1]$.
\end{Lemma}
The proof in \cite{nietic} uses Lemma \ref{L:crit_covarianceII} with a small (and natural) twist: instead of the  original indexing/ordering, for each $N$, $\bm$, to each $k\in \{1,\ldots,N\}$ one assigns $b(k):=\max_{j=1}^d \{a_{k+j}: m_j \neq 0\}$, and then estimates 
$|E(Y_k(\bm) \overline{Y_l}(\bm)+ Y_l(\bm) \overline{Y_k}(\bm))|$ whenever  $\min\{|b(k)-b(l)|, b(k),b(l)\}\geq \sqrt{N}$. 
The rest of the argument is exactly as described above (without the probabilistic rescaling by $a-1$). 

Further qualitative and quantitative improvements and extensions, to be recalled soon, were done in the late 1980s by Niederreiter and Tichy \cite{nietic_real}, Tichy \cite{tichy87},
Drmota, Tichy and Winkler \cite{drmticwin},
and Goldstern \cite{goldstern87}
(see also \cite{drmota_tichy}, Section 1.6).
And yet, interesting  non-trivial open problems remain, as indicated in the next section.

\section{Concluding remarks with open problems}
\label{S:conclOP}
The main theorem of \cite{comed} also yields novel generators of (stochastic a.s.) completely equidistributed sequences in the non-linear setting, of which the prototype is $x_k(t)= (t (\log{t})^{r_k})^{w_k}$, $k\geq 1$, where $(r_k)_k$ (resp.~$(w_k)_k$) are sequences of positive (resp.~natural) numbers satisfying certain hypotheses. 
This may not be so interesting (even though there is no obvious reduction of such sequences to exponential sequences of  Section \ref{S:discNT}) in view of the main theorem of Niederreiter and Tichy \cite{nietic,nietic_real}. However, the approach to (complete) equidistribution, implicit in \cite{comed} and made explicit in Section \ref{S:equi_prob}, is interesting.
In particular, it leads to the realization that \cite{nietic} is based on probability arguments.

\medskip
A careful reader must have noticed that  (\ref{Edelv_SLLN}) is only a sufficiency condition for the SLLN (\ref{Econcl}). 
Davenport, Erd\"os and LeVeque also state (in \cite{delv63}, main/only Theorem) ``On the contrary...", however these counter-examples do not either imply nor disprove the necessity of assumption
$$
\sum_n \frac{1}{n} E\left(\frac{|S_n|}{n}\right)^2 <\infty,
$$
for the SLLN of the corresponding Weyl sums. 
Indeed, it could be that $x_k$ are purely constant (random variables),  and such that
 for some $\bm\in \bbZ^d \setminus {\bf 0}$ we have that 
$\frac{E|S_n(\bm)|^2}{n^2}=\frac{|S_n(\bm)|^2}{n^2}$ is both
$\Omega(1/{\log{n}})$ and $o(1)$ (Korobov \cite{kor60} guarantees that such examples exist), so that the SLLN happens even though the series in the criterion diverges. 
However, here we would trivially have, for the same $d$ and $\bm\in \bbZ^d  \setminus \{\bf 0\}$, that
\begin{equation}
\label{Edelv_SLLNv}
\sum_n \frac{1}{n} {\rm var}\!\left(\frac{|S_n(\bm)|}{n}\right) <\infty, \mbox{ and }  \lim_n \frac{E|S_n(\bm)|^2}{n^2} =0. 
\end{equation}
While the above example is arguably contrived, it already illustrates the following straight-forward consequence of (\ref{Edelv_SLLN}) (or more precisely, of its generalization \cite{lyo88}, Theorem 1):
 (\ref{Edelv_SLLNv}) is a (strictly) stronger sufficiency criterion for 
the SLLN of the corresponding (indexed by $\bm$) Weyl sums, than that given in \cite{delv63} (see (\ref{Edelv_SLLN})). Is it necessary? If not, is the whole $d$-dimensional collection of them (meaning that (\ref{Edelv_SLLNv}) is true  for all $\bm\in \bbZ^d  \setminus \{\bf 0\}$) necessary for $d$-multiple equidistribution of the corresponding $\bbeta$? If not, is the complete collection of them
(meaning that (\ref{Edelv_SLLNv}) is true for all $\bm\in \bbZ^d  \setminus \{\bf 0\}$ and all $d\in \bbN$) necessary for complete equidistribution of the corresponding  $\bbeta$? If not, are there natural additional (easy to check) hypotheses on $(x_k)_k$ under which 
these families of criteria become necessary (they are always sufficient)?

\medskip
As pointed out in Section \ref{S:discNT}, we know from \cite{nietic} that Koksma's numbers are of R4 type and that  the same is true, as shown in \cite{nietic_real}, for a large class of other exponentially generated sequences. 
A trivial example of a c.e.s.~which is also R4 (and even R6 in Knuth's terminology) is the infinite statistical sample $\Gamma$, but as explained in Section \ref{S:why_not}, this is likely not considered sufficiently explicit from a non-probabilist's perspective.
The examples from Section \ref{S:Twonew} must have been known to the modern analytic number theory community, in particular as examples of so-called {\em lacunary sequences}. We could not find them on the list of c.e.s.~in any the standard references, however (\ref{ex_novi2}) was studied by Korobov \cite{kor50} in connection to an explicit example of a simply equidistributed sequence. They seem amenable to analysis, and natural for continuing the investigation in Knuth's framework, who anticipated the result of \cite{nietic}, but at the time knew only the work of Franklin and Koksma about Koksma's numbers (and variations thereof) \cite{koksma,F2}.
Are linear c.e.s.~from Section \ref{S:Twonew} also of R4 type?
Is it possible that any sufficiently random c.e.s.~(e.g.~${\rm var}(\beta_k\bbJ_{[c,d]}) >0$ for all $k$ and and all $0\leq c<d\leq 1$) is of R4 type?
If not, is there a natural and easy to check characterization of when the complete equidistribution (type R1)
implies R4?
By the way, note that if an infinite sequence of random variables $(\beta_k)_k$ on $([0,1],{\mathcal B},\mes)$ is both completely equidistributed in $[0,1]$ (that is, of type R1) and {\em exchangeable} (see e.g.~\cite{durrett}, Example 5.6.4), then it must be equal in law to $\Gamma$. 
Is there a natural exchangeability-like property (clearly stronger than R4 type), though weaker than exchangeability, that would jointly with  c.e.s.~still imply that the sequence $(\beta_k)_k$ has the law of $\Gamma$?
Is there a natural condition (stronger than c.e.s.~but weaker than i.i.d.) that would imply,  jointly with R4 type, that the sequence $(\beta_k)_k$ has the law of $\Gamma$?

\medskip
Once equidistribution is established, it is natural to ask about the rate of convergence in (\ref{EeqWeyl}), or in analogous multi-dimensional expressions that correspond to Definitions \ref{def1}-\ref{def2}.
The discrepancy sequences (as recalled in Section \ref{S:equiv_form}) seem to be the principal object for this analysis, already in the focus of classical studies by the founders of analytic number theory (see \cite{kuni}, Chapter 2 and \cite{drmota_tichy} Sections 1.1-1.2). 
From the stochastic point of view, it seems more interesting to study discrepancy as a random variable. In particular, the well-known
{\em Erd\H{o}s-Tur\'an-Koksma inequality} (see e.g.~\cite{drmota_tichy}, Theorem 1.21) can be restated in the present setting as follows: if $H$ is a positive integer 
and, for $\bm \in \bbZ^d$, we let $r(\bm):=\prod_{i=1}^d\max\{1,|m_i|\}$,
then it is true ($\bt$-by-$\bt$) that
\begin{eqnarray}
\nonumber 
D_{N}^{*,d}&\leq&
\left(\frac{3}{2}\right)^d
\left(
\frac{2}{H+1}+
\sum_{0<\|\bm\|_{\infty}\leq H}\frac{1}{r(\bm)}
\left|
\frac{1}{N}
\sum_{k=1}^{N} Y_k(\bm)
\right|
\right)\!.
\end{eqnarray}
Substantial work by the ``Tichy school'' has been done on multi-dimensional discrepancy estimation
of Koksma's numbers and variations.
In particular, Tichy \cite{tichy87} obtains, in the context of \cite{nietic_real} (more precisely, $x_k(t)=t^{a_k}$, $t>1$, for some $a_k\in {\mathbb R}$, and the minimal distance between different powers is bounded below by $\delta>0$), bounds on $D_{N}^{*,d}$ of the form $O(N^{-1/2+\eta})$ for any $\eta>0$, even if the multiplicity $d$ is not fixed but diverges as $\log{N}$ raised to a small power; and 
 Goldstern \cite{goldstern87} extends these to the setting where the replication between the exponents $a_n$ is possible but infrequent, and the minimal distance between the exponents may slowly converge to $0$ (see also notes on the literature in \cite{drmota_tichy}, Section 1.6). 

Recall that already in the setting of linear generators,  
the random variables $Y_k(\bm)$ are not mutually independent, however they do have other nice properties. 
For concreteness, one could take the two examples of completely equidistributed sequences from Section \ref{S:Twonew}.
Is there a central limit theorem, or another type of concentration result that would apply, and give good estimates on $d$-dimensional (star) discrepancy of these sequences, with or without modifications related to R4 type of randomness?
When is 
$
\sqrt{N/\log{\log{N}}} D_N^{*,d}
$
a tight sequence of random variables?
The above LIL-type result for the simple ($1$-multiple) discrepancy sequences $(D_N^{*,1})_N$ is well-known  (see e.g.~\cite{drmota_tichy} Section 1.6.2 or \cite{aisber_review}) in the context of lacunary sequences, and in particular for  (\ref{ex_novi1},\ref{ex_novi2}), even if permuted (without possible deletion), as shown by Aistleitner et al.~\cite{aisbertic}. But as soon as one starts increasing the multiplicity $d$, no specific study of the corresponding LIL seem to exist.

\medskip
In view of (\ref{Eequiv_form_Gamma}) and Remark 1(d), most probabilists and statisticians would likely ask if any of the following is true for any (or all) of the c.e.s.~discussed here: given any $d\geq 1$ and any $f:\bbR^d \to \bbR$ a nice enough function, define $f(\bbeta_k^\infty):=f(\bbeta_k^d)$ and $f(\Gamma):=f(U_1,U_2,\ldots,U_d)$. Are the sequences of random variables (recall (\ref{ex1.6mult}))
\begin{equation}
\label{ECLT1}
\sqrt{N}\left(\frac{1}{N} \sum_{k=0}^{N-1} f(\bbeta_{kh+1}^\infty) - {\bf{E}}\, (f(\Gamma))\right), \ N\geq 1, \mbox{ where } h\geq d,
\end{equation}
and/or
\begin{equation}
\label{ECLT2}
\sqrt{N}\left(\frac{1}{N} \sum_{k=1}^{N} f(\bbeta_{k}^\infty(t)) - {\bf{E}}\, (f(\Gamma))\right), N \geq 1,
\end{equation}
tight? Do they converge in law to a centered Gaussian random variable? 
These questions are similar in spirit to those of the preceding paragraph, but not quite the same (see \cite{aisber_funcCLT} for a functional CLT for 1-dimensional discrepancies). 
In particular, the well-known {\em Koksma} and {\em Koksma-Hlawka} inequalities (see e.g.~\cite{kuni}, Ch.~2, Section 5 or \cite{drmota_tichy}, Theorem 1.14) serve to give universal bounds on the error in (MC) numerical integration in terms of the discrepancy of the sequence and (a multi-dimensional extension of) the bounded variation of the integrand. Though essential in various applications, they are too crude for studying weak convergence properties (\ref{ECLT1},\ref{ECLT2}). 
In the one-dimensional setting, again a substantial progress has been made for lacunary sequences, starting from the classical work of Salem and Zygmund \cite{salem_zygmund}, and ending in the recent studies by Aistleitner et al.~\cite{aisber_CLT,aisbertic} (see also \cite{aistleitner_MG,aisber_review}).
Satisfying the  CLT analogues (\ref{ECLT1},\ref{ECLT2}) (with permutations and possible deletions permitted) would undoubtedly be a strong ``evidence of randomness''.  How is it (or is it) related to the strongest Knuth's  \cite{knu69} R6 type of pseudo-randomness? 

\medskip
{\bf Acknowledgments}  V.L.~thanks Chris Burdzy for very useful pointers to the literature, and several constructive comments that reset this research project on the right track, and Martin Goldstern and Pete L.~Clark for additional help with the literature. We are both grateful to Christoph Aistleitner for detailed reading of the original preprint, and a number of valuable comments and suggestions. 
The first version of this paper was written while V.L. was at the Universit\'e Paris-Sud XI.

\end{document}